\documentclass[12pt,a4paper,twoside]{article}
\usepackage{amsmath,amsfonts,amssymb,a4,enumerate,epsfig,array,theorem,psfrag}

\newtheorem{theo}{Theorem}
\newtheorem{prop}{Proposition}[section]
\newtheorem{coro}[prop]{Corollary}

\newtheorem{lemma}[prop]{Lemma}

\newcommand{\ZZ}{{\mathbb{Z}}}
\newcommand{\QQ}{{\mathbb{Q}}}
\newcommand{\RR}{{\mathbb{R}}}

\newcommand{\tD}{{\tilde D}}
\newcommand{\Gg}{{\cal{G}}}
\newcommand{\cG}{{\cal{G}}}

\newcommand{\cT}{{\cal{T}}}
\newcommand{\cD}{{\cal{D}}}
\newcommand{\cR}{{\cal{R}}}

\newcommand\rank{\operatorname{rank}}
\newcommand\step{\operatorname{step}}
\newcommand\diag{\operatorname{diag}}

\newcommand{\nobf}{\noindent\bf}

\renewcommand\emptyset{\varnothing}
\def\qed{\unskip\nobreak\hfil\penalty50\hskip1.75em\null\nobreak\hfil
$\blacksquare$ {\parfillskip=0pt \finalhyphendemerits=0 \par}\goodbreak}

\begin{document}
\title{Cut-and-paste of quadriculated disks and \\
arithmetic properties of the adjacency matrix}
\author{Nicolau C. Saldanha and Carlos Tomei}
\maketitle

\begin{abstract}
We define \textit{cut-and-paste},
a construction which, given a quadriculated disk
obtains a disjoint union of quadriculated disks of smaller total area.
We provide two examples of the use of this procedure as a recursive step.
Tilings of a disk $\Delta$ receive a \textit{parity}:
we construct a perfect or near-perfect matching
of tilings of opposite parities.
Let $B_\Delta$ be the black-to-white adjacency matrix:
we factor $B_\Delta = L\tilde DU$,
where $L$ and $U$ are lower and upper triangular matrices,
$\tilde D$ is obtained from a larger identity matrix
by removing rows and columns and
all entries of $L$, $\tilde D$ and $U$ are
equal to $0$, $1$ or $-1$.
\end{abstract}

\section{Introduction}

\begin{figure}[ht]
\begin{center}
\epsfig{height=3cm,file=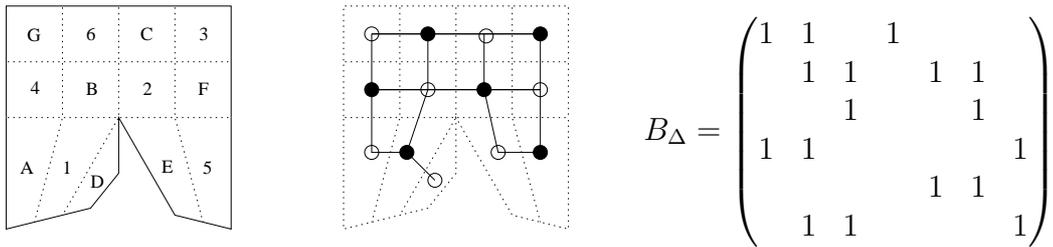}
\qquad
\raise 12mm \hbox{\(
B_\Delta = \begin{pmatrix}
1 & 1 &  & 1 &  &  &  \\
 & 1 & 1 &  & 1 & 1 &  \\
 &  & 1 &  &  & 1 &  \\
1 & 1 &  &  &  &  & 1 \\
 &  &  &  & 1 & 1 &  \\
 & 1 & 1 &  &  &  & 1
\end{pmatrix}
\)}
\end{center}
\label{fig:disk0}
\caption{A quadriculated disk, its dual graph and its black-to-white matrix}
\end{figure}

\footnote{2000 {\em Mathematics Subject Classification}.
Primary 05B45, 05C70; Secondary 05B20, 05C50.
{\em Keywords and phrases} Quadriculated disk,
matchings, tilings by dominoes, dimers.}
\footnote{The authors gratefully acknowledge the support of CNPq, Faperj.
The first author thanks the kind hospitality of The Ohio State University
during part of the time when this paper was written.}
In this paper, a \textit{square} is a topological disk
with four privileged boundary points, the \textit{vertices};
the boundary of the square consists of four \textit{edges}.
A \textit{quadriculated disk} $\Delta$ is a closed topological disk
formed by the juxtaposition along edges of finitely many squares
such that interior vertices belong to precisely four squares:
it may be considered
as a closed subset of the plane $\RR^2$ tiled by quadrilaterals.
A simple example is the $n \times m$ rectangle divided into unit squares,
another is shown in Figure~\ref{fig:disk0}.

Given $\Delta$, we define the planar \textit{dual graph} $\Gg_\Delta$:
vertices of $\Gg_\Delta$ correspond to squares in $\Delta$ and
two vertices of $\Gg_\Delta$ are adjacent if their corresponding
squares share an edge.
Quadriculated disks are {\it bi-colored}:
the squares are black and white in a way that
squares with a common edge have opposite colors
(equivalently, $\Gg_\Delta$ is bipartite).
Label the black (resp. white) squares of a quadriculated disk $\Delta$
by $1, 2, \ldots, b$ (resp. $1, 2, \ldots, w$).
The $b \times w$ {\it black-to-white} (adjacency) matrix $B_\Delta$
has $(i,j)$ entry $b_{ij} = 1$ if the $i$-th black and $j$-th white squares
share an edge and $b_{ij} = 0$ otherwise.
Figure \ref{fig:disk0} is an example of black-to-white matrix;
black and white squares are labeled by numbers and letters, respectively.
Throughout the paper, blank matrix entries equal $0$.
For a labeling in which black vertices come first,
the adjacency matrix  of $\Gg_\Delta$ is
\[ \begin{pmatrix} 0 & B_\Delta \\ B_\Delta^T & 0 \end{pmatrix}. \]
The following result \cite{DT} indicates an unexpected spectral rigidity
of $B_\Delta$.

\begin{theo}
\label{theo:DT}
Let $\Delta$ be a quadriculated disk with $b=w$ and black-to-white
matrix $B_\Delta$. Then $\det(B_\Delta)$ equals $0$, $1$ or $-1$.
\end{theo}

This result admits a combinatorial interpretation.
A {\it domino tiling} $\tau$ of $\Delta$ is a decomposition of $\Delta$
as a union of \textit{dominos} (i.e., $2 \times 1$ rectangles)
with disjoint interior.
Let $\cT_\Delta$ be the set of domino tilings of $\Delta$.
There is a natural parity function on $\cT_\Delta$ (see Section 4) and
the determinant $\det(B_\Delta)$ counts tilings with a sign given by parity.
The theorem above thus says that there exists a \textit{quasi-perfect matching}
in $\cT_\Delta$, i.e., a correspondence between even and odd tilings
leaving out at most one element of $\cT_\Delta$, the \textit{loner}.
We provide a new, (quasi-) bijective proof of Theorem \ref{theo:DT}
by constructing a quasi-perfect matching in the bipartite set $\cT_\Delta$.

We extend Theorem \ref{theo:DT} in a different, more algebraic, direction.
A rectangular matrix $\tD$ is a \textit{defective identity}
if it can be obtained from the identity matrix
by adding rows and columns of zeros.
For a $n \times m$ matrix $A$,
an $L\tD U$ decomposition of $A$ is a factorization $A = L\tD U$
where $L$ (resp. $U$) is $n \times n$ (resp. $m \times m$)
lower (resp. upper) invertible and $\tD$ is a defective identity.


\begin{theo}
\label{theo:LDU}
Let $\Delta$ be a quadriculated disk with at least two squares.
For an appropriate labeling of its squares,
the black-to-white matrix $B_\Delta$
admits an $L\tD U$ decomposition
whose factors have all entries equal to $0$, $1$ or $-1$.
\end{theo}

Thus, for example, the matrix $B_\Delta$ in Figure \ref{fig:disk0}
admits the decomposition
\[
\begin{pmatrix}
1 & & & & & \\
& 1 & & & & \\
& & 1 & & & \\
1 & & & -1 & & \\
& & & & 1 & \\
& 1 & & & -1 & 1 
\end{pmatrix}
\begin{pmatrix}
1 & & & & & & \\
& 1 & & & & & \\
& & 1 & & & & \\
& & & 1 & & & \\
& & & & 1 & 0 & \\
& & & & & & 1
\end{pmatrix}
\begin{pmatrix}
1 & 1 & & 1 & & & \\
& 1 & 1 & & 1 & 1 & \\
& & 1 & & & 1 & \\
& & & 1 & & & -1 \\
& & & & 1 & 1 & \\
& & & & & 1 & \\
& & & & & & 1
\end{pmatrix}.
\]

Both the construction of the quasi-perfect matching and the proof
of Theorem \ref{theo:LDU} use \textit{cut-and-paste},
a recursive operation on quadriculated disks.
A quadriculated disk $\Delta$ is \textit{cut} along \textit{diagonals}
and \textit{pasted} to obtain
a disjoint union of smaller disks $\Delta'_1, \ldots, \Delta'_d$,
often with $d = 1$ (Lemma \ref{lemma:deltalinha}).
Every nontrivial quadriculated disk admits cut-and-paste
(Proposition \ref{prop:gooddiag}).


The proof of Theorem \ref{theo:LDU} relies on
a procedure to convert $L\tilde DU$ decompositions of
$B_{\Delta'_1}, \ldots, B_{\Delta'_d}$
into a similar decomposition of $B_\Delta$ (Lemma \ref{lemma:LDUstep}).
The proof yields a fast algorithm to obtain the appropriate labeling
of vertices, the matrices in the factorization,
$\det(B_\Delta)$ and $\rank(B_\Delta)$.

In Section 2 we present the facts about \textit{diagonals}
of quadriculated disks which will be used in
Section 3 to describe \textit{cut-and-paste}.
In Section \ref{sect:bijective} we construct the quasi-perfect matching.
The inductive step in the proof of Theorem \ref{theo:LDU},
the algebraic counterpart of cut-and-paste,
is the main topic of Section \ref{sect:decomp}.
Finally, in Section \ref{sect:boards},
we study \textit{boards},
quadriculated disks which are subsets of the quadriculated plane
$\ZZ^2 \subset \RR^2$;
Theorem \ref{theo:board} states that cut-and-paste can be performed within
this smaller class.

Counting tilings with sign given by parity (as in Theorem \ref{theo:DT})
corresponds to the case $q=-1$ of 
the $q$-counting of domino tilings with respect
to \textit{height} or \textit{volume}
as in \cite{T}, \cite{EKLP} and \cite{STCR}. 
In a similar vein, \cite{ST} extends Theorem \ref{theo:DT}
to quadriculated annuli by introducing a polynomial
which counts tilings with respect to yet another integral
parameter, the \textit{flux}.
It is not clear whether the cut-and-paste procedure
can be extended to take such parameters into account.


\section{Diagonals}
\label{sect:diag}

\begin{figure}[ht]
\begin{center}
\psfrag{d1}{$\delta_1$}
\psfrag{d2}{$\delta_2$}
\psfrag{d3}{$\delta_3$}
\psfrag{d4}{$\delta_4$}
\psfrag{d5}{$\delta_5$}
\psfrag{d6}{$\delta_6$}
\psfrag{v0}{$v_0$}
\psfrag{v1}{$v_1$}
\psfrag{v2}{$v_2$}
\psfrag{s12}{$s_{\frac{1}{2}}$}
\psfrag{s32}{$s_{\frac{3}{2}}$}
\epsfig{height=35mm,file=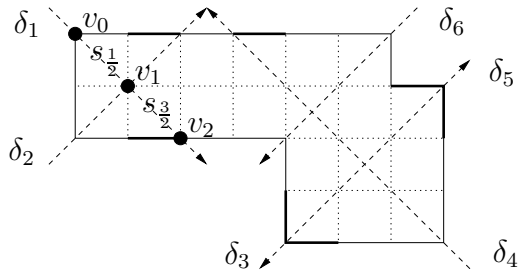}
\end{center}
\caption{A quadriculated disk and its six diagonals}
\label{fig:gooddiag}
\end{figure}

A \textit{corner} of a quadriculated disk $\Delta$
is a boundary point which is a vertex of a single square.
A \textit{pre-diagonal} of length $k > 0$ of $\Delta$ is
a sequence of vertices $v_0v_1\ldots v_k$ such that
\begin{enumerate}[(i)]
\item{$v_0$ is a corner, $v_1, v_2, \ldots, v_{k-1}$ are interior vertices;}
\item{consecutive vertices $v_i$ and $v_{i+1}$, $i = 0, \ldots, k-1$,
are opposite vertices of a square $s_{i+1/2}$;}
\item{consecutive squares $s_{i-1/2}$ and $s_{i+1/2}$,
$i = 1, \ldots, k-1$, have a single vertex in common (which is $v_i$);}
\item{the vertices $v_i$ and the squares $s_{i+1/2}$, $i = 0, \ldots, k-1$,
are distinct.}
\end{enumerate}
A \textit{diagonal} is a maximal pre-diagonal (under inclusion).
More geometrically, we may think of a diagonal as a line
$\ell = \ell(v_0,s_{1/2},v_1,\ldots ,s_{k-1/2},v_k)$
connecting $v_0$, the center of $s_{1/2}$,
$v_1$, the center of $s_{3/2}$ and so on up to $v_k$.
The squares $s_{1/2}, \ldots, s_{k-1/2}$
are the \textit{squares of the diagonal}.
Usually, the vertices $s_{1/2}, \ldots, s_{k-1/2}$
form a cut set of the dual graph $\cG_\Delta$.
Diagonals, being sequences of vertices, are naturally oriented.
Figure \ref{fig:gooddiag} shows examples of diagonals;
vertices and squares of $\delta_1$ are indicated.

\begin{prop}
\label{prop:simplediags}
Given a corner $v_0$ of $\Delta$ there is a unique diagonal
starting at $v_0$.
Furthermore, all diagonals end at boundary points.
\end{prop}

\begin{figure}[ht]
\begin{center}
\psfrag{v0}{$v_0$}
\psfrag{vk}{$v_k$}
\psfrag{vi}{$v_i$}
\psfrag{vi=vk}{$v_i = v_k$}
\psfrag{(a)}{(a)}
\psfrag{(b)}{(b)}
\epsfig{height=30mm,file=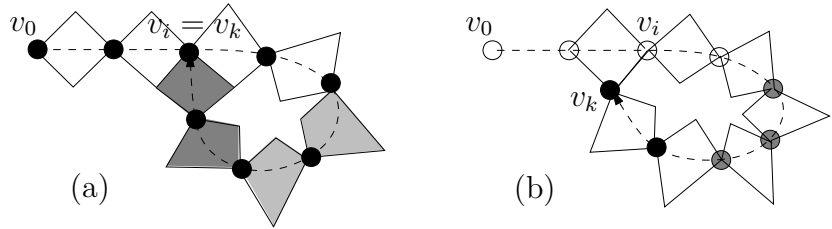}
\end{center}
\caption{Impossible diagonals}
\label{fig:prediag2}
\end{figure}

{\nobf Proof:}
In principle, there are three types of diagonals:
the vertex $v_k$ may coincide with some $v_i$, $i < k$
(Figure \ref{fig:prediag2}, (a)),
the square $s_{k+1/2}$ may coincide with some $s_{i+1/2}$, $i <k$
(Figure \ref{fig:prediag2}, (b))
or $v_k$ may be a boundary vertex of $\Delta$.
Existence and uniqueness of a diagonal $\delta$ starting at the corner $v_0$
follows from finiteness.
The reader may check that self-intersection would
happen at right angles, as in the figure.
Bicoloring of squares and vertices of $\Delta$,
as in Figure \ref{fig:prediag2}, yields a contradiction in either case.
\qed

Let $\delta$ be a diagonal of a quadriculated disk $\Delta$
associated to the line $\ell = \ell(v_0, s_{1/2}, \ldots, v_k) \subset \Delta$.
Given a vertex $v$ of $\Delta \smallsetminus \ell$,
draw a smooth curve $\gamma: [0,1] \to \Delta$,
$\gamma(0) = v$, $\gamma(1) \in \ell$,
$\gamma(t) \in \Delta \smallsetminus \ell$ for $t < 1$
and $\gamma'(1)$ transversal to $\ell$.
We say that $v$ is to the \textit{left} (resp. \textit{right})
of $\delta$ if $\det(v_1 - v_0, \gamma'(1))$ is negative (resp. positive).
The existence of the curve $\gamma$ follows from 
the fact that $\Delta$ is path-connected.
A vertex $v$ is not simultaneously to the left and right of $\delta$:
indeed, if $\gamma_l, \gamma_r: [0,1] \to \Delta$
satisfy the hypothesis above and
$\det(v_1-v_0,\gamma_l'(1)) < 0 < \det(v_1-v_0,\gamma_r'(1)$
then juxtaposition of $\gamma_l$ and time-reversal of $\gamma_r$
obtains a loop which crosses $\ell$ exactly once, a contradiction.

In the next section we will use diagonals to cut-and-paste.
Not all diagonals are suitable for this construction.
Call the two edges of $s_{k-1/2}$ ending at $v_k$ 
\textit{terminal} edges.
A diagonal $v_0\ldots v_k$ is a \textit{good diagonal}
if at least one terminal edge is contained 
in the boundary of $\Delta$.
In Figure \ref{fig:gooddiag}, $\delta_6$ is the only bad diagonal.
A square has four diagonals, all good.


To prove the existence of good diagonals,
we use a quadriculated version of the Gauss-Bonnet theorem.
Let $V$ be the number of vertices of $\Delta$ and write
$V = V_I + V_1 + V_2 + \cdots + V_r$
where $V_I$ counts interior vertices and
$V_r$ is the number of boundary vertices
belonging to exactly $r$ squares.
Notice that $V_1$ is the number of corners of $\Delta$.

\begin{lemma}
\label{lemma:gaussbonnet}
$V_1 - V_3 - 2V_4 - \cdots - (r-2)V_r = 4$.
\end{lemma}

{\nobf Proof:}
Let $E$ and $F$ be the number of edges and faces
(i.e., squares) of $\Delta$.
Write $E = E_I + E_B$, where $E_I$ (resp. $E_B$)
counts interior (resp. boundary) edges.
Clearly, $4F = 2E_I + E_B = 2E - E_B$
and therefore $4E = 8F + 2E_B$.
Also, $4F = 4 V_I + V_1 + 2 V_2 + \cdots + r V_r =
4V - (3 V_1 + 2 V_2 + \cdots + (4-r) V_r)$
and $4V = 4F + (3V_1 + 2V_2 + \cdots + (4-r) V_r)$.
By Euler, $4V - 4E + 4F = 4$.
Substituting the above formulas and using
$E_B = V_1 + V_2 + \cdots + V_r$ we have
the desired identity.
\qed

\begin{prop}
\label{prop:gooddiag}
Any quadriculated disk $\Delta$ admits at least four good diagonals.
\end{prop}

{\nobf Proof:}
Each vertex counted in $V_1$ is a starting corner for a diagonal:
we have to prove that at least four of these $V_1$ diagonals are good.
Each vertex counted in $V_3$, for example, is the endpoint
of at most three diagonals of which only one is declared bad.
More generally, we have at most $V_1 - 4 = V_3 + 2V_4 + \cdots + (r-2)V_r$
bad ends and we are done.
\qed

\section{Geometric cut-and-paste}
\label{sect:cutnpaste}

We are ready to perform {\it cut-and-paste} along a good diagonal.
A good diagonal $v_0\ldots v_k$ is \textit{balanced}
if exactly one terminal edge is contained in the boundary of $\Delta$.
Diagonals $\delta_1$, $\delta_2$ and $\delta_4$ in Figure
\ref{fig:gooddiag} are balanced;
$\delta_3$ and $\delta_5$ are unbalanced.

\begin{figure}[ht]
\begin{center}
\psfrag{del}{$\delta$}
\psfrag{par}{$\partial_\delta$}
\epsfig{height=33mm,file=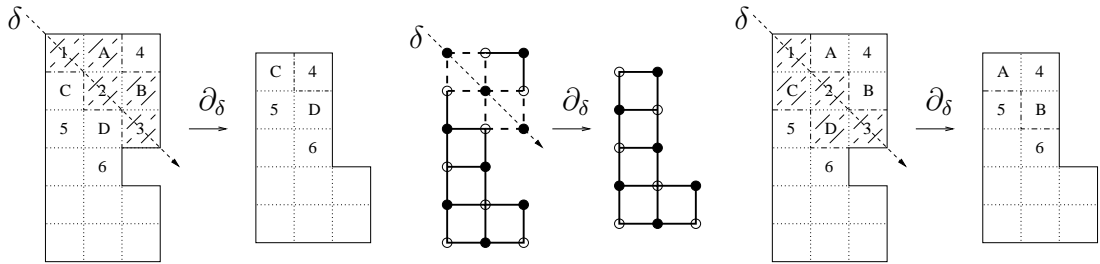}
\end{center}
\caption{Cut-and-paste along the unbalanced diagonal $\delta$ of length $k = 3$}
\label{fig:cutnpaste0}
\end{figure}

In Figure \ref{fig:cutnpaste0} we illustrate the cut-and-paste procedure
$\partial_\delta$ on a quadriculated disk $\Delta$
and its dual graph $\cG_\Delta$, where $\delta$ is an unbalanced diagonal.
The operation removes the shaded squares
and identifies edges to obtain a new quadriculated disk
$\Delta' = \partial_\delta(\Delta)$.
Another choice of shaded squares for the same good diagonal $\delta$
is indicated in the right and obtains the same quadriculated disk $\Delta'$.
In the left (resp. right), squares $C$ and $D$ (resp. $A$ and $B$)
take over the space vacated by $A$ and $B$ (resp. $C$ and $D$).

\begin{figure}[ht]
\begin{center}
\psfrag{del}{$\delta$}
\psfrag{par}{$\partial_\delta$}
\epsfig{height=28mm,file=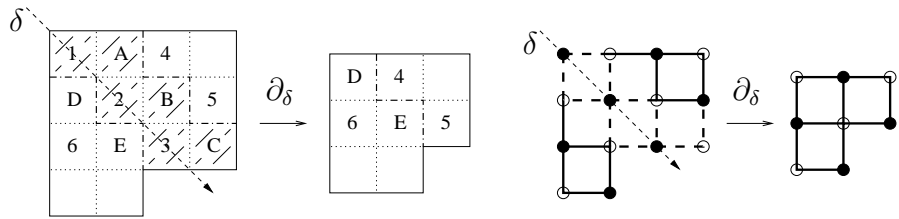}
\end{center}
\caption{Cut-and-paste along a balanced diagonal, $k = 3$}
\label{fig:cutnpaste1}
\end{figure}

The balanced case shown in Figure \ref{fig:cutnpaste1} is a little different.
It turns out that a similar construction with another choice of zig-zag
is not appropriate for future purposes.

In the dual graph $\cG_\Delta$, cut-and-paste removes the cut set
of vertices (of the graph) associated with squares of the diagonal $\delta$
and identifies vertices on both sides:
vertices left without partners at the end of an unbalanced diagonal
are also deleted.
This point of view is more symmetric and does not require
the specification of zig-zags.

Notice that the extreme vertex $v_k$ of a balanced diagonal
may belong to more than two squares, as in Figure \ref{fig:disk0}.
This is innocuous, as we shall see.

\begin{figure}[ht]
\begin{center}
\epsfig{height=35mm,file=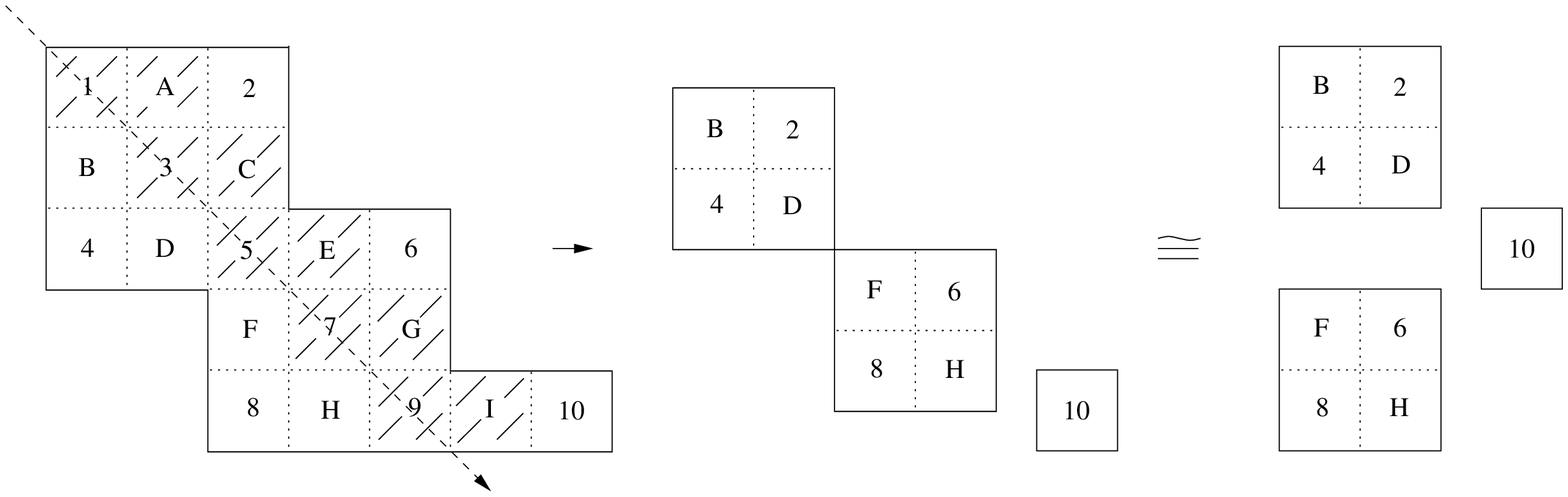}
\end{center}
\caption{Cut-and-paste may produce a disjoint union of disks}
\label{fig:cutnpaste2}
\end{figure}

Cut-and-paste allows for recursive proofs and constructions
in the class of finite disjoint unions of quadriculated disks.
As we shall prove in Lemma \ref{lemma:deltalinha},
given a quadriculated square $\Delta$ and a good diagonal $\delta$,
cut-and-paste obtains a quadriculated region $\tilde\Delta'$
which consists of quadriculated disks
$\Delta'_1, \ldots, \Delta'_d$, possibly joined by points.
The process of passing from $\tilde\Delta'$ to
$\Delta' = \Delta'_1 \sqcup \cdots \sqcup \Delta'_d$
is called \textit{detaching}.
Clearly, $\Delta'$ has fewer squares than $\Delta$.
In the two previous examples, $d = 1$; in Figure \ref{fig:cutnpaste2}, $d = 3$.

In a somewhat degenerate case, $\Delta' = \emptyset$ if and only if
$\Delta$ consists of one or two squares.
Also, if $\delta$ is an unbalanced diagonal of length $k = 1$,
the quadriculated disk $\Delta'$ is obtained from $\Delta$
by deleting two squares.

\begin{figure}[ht]
\begin{center}
\psfrag{del}{$\delta$}
\psfrag{par}{$\partial_\delta$}
\psfrag{void}{$\emptyset$}
\epsfig{height=15mm,file=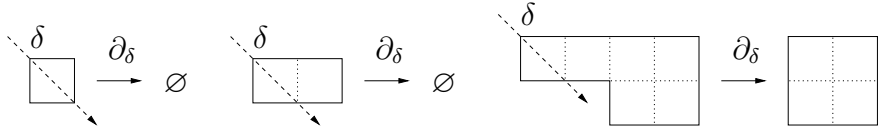}
\end{center}
\caption{Cut-and-paste in extreme situations}
\label{fig:void}
\end{figure}


\begin{lemma}
\label{lemma:deltalinha}
Let $\Delta$ be a quadriculated disk
with a good diagonal $\delta$ of length $k > 1$.
Let $\Delta'$ be obtained from $\Delta$ by cut-and-paste
along $\delta$ (and detaching):
$\Delta'$ is a disjoint union of quadriculated disks.
\end{lemma}

We use a notation for vertices and squares near a good diagonal.
Squares immediately to the left (resp. right) of the diagonal
are labelled $s^l_1, s^l_2, \ldots$ (resp. $s^r_1, s^r_2, \ldots$).
Similarly, vertices to the left (resp. right)
are labelled $v^l_{1/2}, v^l_{3/2}, \ldots$
(resp. $v^r_{1/2}, v^r_{3/2}, \ldots$).
Thus, in Figure \ref{fig:cutnpaste0},
$s^l_1 = A$, $s^l_2 = B$, $s^r_1 = C$, $s^r_2 = D$;
in Figure \ref{fig:cutnpaste1},
$s^l_1 = A$, $s^l_2 = B$, $s^l_3 = C$, $s^r_1 = D$, $s^r_2 = E$.
The squares deleted in the cut-and-paste construction
(dashed in the figures)
are $s_{1/2}, s^x_1, s_{3/2}, \ldots, s^x_{k-1}, s_{k-1/2}$
and, in the balanced case, $s^x_k$; here $x = l$ or $x = r$.
Let $\Delta^r$ (resp. $\Delta^l$) be the closed regions
to the right (resp. left) of the deleted squares.
Attach $\Delta^l$ to $\Delta^r$ by identifying edges
in order to obtain a quadriculated region $\tilde\Delta'$.



\begin{figure}[ht]
\begin{center}
\psfrag{s12}{$s_{\frac{1}{2}}$}
\psfrag{s32}{$s_{\frac{3}{2}}$}
\psfrag{s52}{$s_{\frac{5}{2}}$}
\psfrag{sl1}{$s^l_1$}
\psfrag{sl2}{$s^l_2$}
\psfrag{sl3}{$s^l_3$}
\psfrag{sr1}{$s^r_1$}
\psfrag{sr2}{$s^r_2$}
\psfrag{sr3}{$s^r_3$}
\psfrag{Dl}{$\Delta^l$}
\psfrag{Dr}{$\Delta^r$}
\psfrag{del}{$\delta$}
\psfrag{lr}{$\ell^r$}
\psfrag{(a)}{(a)}
\psfrag{(b)}{(b)}
\epsfig{height=4cm,file=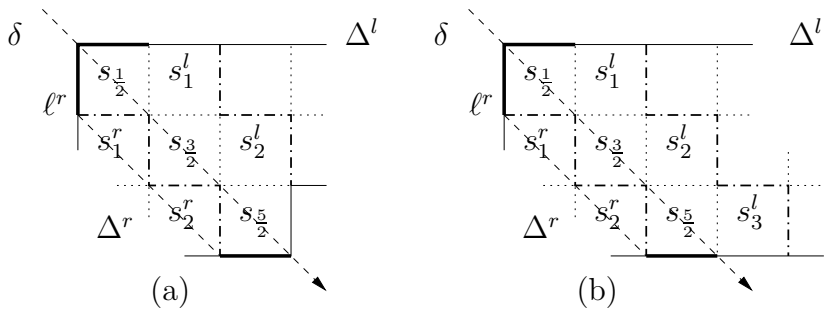}
\end{center}
\caption{Notation for cut-and-paste; unbalanced and balanced cases}
\label{fig:indices}
\end{figure}

{\nobf Proof:}
Assume without loss that cut-and-paste along $\delta$
deletes the squares $s_{1/2}, s^l_1, \ldots, s^l_{k-1}, s_{k-1/2}$ and,
if $\delta$ is unbalanced, $s^l_k$.
We claim that $\Delta^r$ is non-empty, path-connected and simply connected.
Indeed, the squares $s^r_1, \ldots, s^r_{k-1}$ exist
(since $v_1, \ldots, v_{k-1}$ are interior points, $k > 1$).
To show that $\Delta^r$ is path-connected,
it suffices to join by a path in $\Delta^r$
any point $x \in \Delta^r$ to 
the line $\ell^r = (v^r_{1/2},s^r_1, \ldots, s^r_{k-1}, v^r_{k-1/2})$.
Notice that the edges $v_0v^r_{1/2}$ and $v_kv^r_{k-1/2}$
are in the boundary of $\Delta$.
If $x \in \Delta^r$ lies between $\ell$ and $\ell^r$ then $x$ belongs to 
one of the squares $s^r_1, \ldots, s^r_{k-1}$ and the path is easy to construct.
Otherwise, take $\gamma: [0,1] \to \Delta$
as in the definition of left and right of $\delta$ in Section \ref{sect:diag};
$\gamma$ must cross $\ell^r$ and a restriction of $\gamma$
yields the required path.
As to simple connectivity, take a simple closed curve $\alpha$
contained in $\Delta^r$ and therefore in $\Delta$.
By Jordan's Theorem, $\alpha$ encloses a disk $A$.
Since $\Delta$ is simply connected, $A \subset \Delta$.
Also, a path in $A$ from $x \in A$ to $\alpha$
guarantees that $x$ and $\alpha$ are on the same side of $\delta$.

The region $\Delta^l$ may be disconnected or even empty.
On the other hand, the argument above shows that
its connected components are simply connected.
Thus, $\tilde\Delta'$ is obtained by gluing the simply connected pieces
$\Delta^r$ and the components of $\Delta^l$:
we must now study the gluing process.
Let $\zeta^r$ and $\zeta^l$ be the zig-zag lines
$v_{1/2}^r v_1 v_{3/2}^r \ldots v_{k-1} v_{k-1/2}^r$ and
$v_{1/2}^l v_1^{ll} v_{3/2}^l \ldots v_{k-1}^{ll} v_{k-1/2}^l$,
where $v_i^{ll}$ is the left-most vertex of $s_i^l$.
Cut-and-paste obtains $\tilde\Delta'$ by gluing
$\Delta^r$ and $\Delta^l$ along $\zeta^r$ and $\zeta^l$.
Notice that $\zeta^r$ is contained in the boundary of $\Delta^r$.
It is possible, however, that parts of $\zeta^l$
are part of the boundary of $\Delta$ and not in $\Delta^l$.

We claim that, given a connected component $D$ of $\Delta^l$,
its intersection with $\zeta^l$ is either empty or path-connected.
In other words, for any two points $x_0, x_1 \in D \cap \zeta^l$,
the segment $[x_0,x_1] \subset \zeta^l$ between $x_0$ and $x_1$
is contained in $D$.
Indeed, there is a curve $\alpha$ in $D$ joining $x_0$ and $x_1$.
Juxtaposition of $\alpha$ and $[x_0,x_1]$ obtains a closed curve in $\Delta$.
As before, simple connectivity of $\Delta$ implies that the region
surrounded by this closed curve is contained in $\Delta$ and therefore
in $\Delta^l$ and $D$, completing the proof of the claim.

The claims and Seifert-Van Kampen's Theorem (\cite{Massey})
imply that each connected component of $\tilde\Delta'$ is simply connected.
Detaching guarantees that each connected component of $\Delta'$
is a simply connected surface with boundary --- a disk.
\qed


\section{A bijective proof of Theorem \ref{theo:DT}}
\label{sect:bijective}

A nonzero entry $b_{ij}$ of the black-to-white matrix $B_\Delta$
corresponds to a domino contained in $\Delta$:
the indices $i$ and $j$ indicate the black and white squares
in the domino and $b_{ij} \ne 0$ when these two squares are adjacent.
A \textit{domino tiling} of $\Delta$
is a decomposition of $\Delta$ as a union of dominos
with disjoint interiors;
let $\cT_\Delta$ be the set of all domino tilings of $\Delta$.
A nonzero monomial of the black-to-white matrix $B_\Delta$
corresponds to some $\tau \in \cT_\Delta$.
Indeed, the dominos associated with the entries cover $\Delta$
and their interiors are disjoint.
Equivalently, for a labeling of black and white squares
by $\{1, 2, \ldots, b\}$ and $\{1, 2, \ldots, w\}$, 
we may consider a tiling $\tau$ as a function 
$\pi: \{1, 2, \ldots, w\} \to \{1, 2, \ldots, b\}$
with $\pi(j) = i$ if and only if the $i$-th black square
and the $j$-th white square form a domino in $\tau$.
With $b=w$, this provides an identification between $\cT_\Delta$
and a subset of the symmetric group $S_w$.

The above identification endows a tiling with parity (or sign).
Tilings differing by a \textit{flip}
(i.e., by exactly two dominos forming a $2 \times 2$ square)
have opposite parities:
if their corresponding permutations are $\pi_1$ and $\pi_2$
then $\pi_2^{-1}\pi_1$ is a cycle of length $2$,
interchanging the two white squares in the flip.
The combinatorial interpretation of
Theorem \ref{theo:DT} is that the number of even and odd tilings
in $\cT_\Delta$ differ by at most $1$.
In this section we provide a bijective proof of this statement.

More precisely, we present an algorithm that,
given a quadriculated disk $\Delta$, obtains a \textit{quasi-perfect matching}
in $\cT_\Delta$, i.e., a subset $\cT^\ast_\Delta \subseteq \cT_\Delta$
whose complement has at most one element, the \textit{loner},
and an involution $\rho: \cT^\ast_\Delta \to \cT^\ast_\Delta$ 
(i.e., $\rho^2(\tau) = \tau$) inverting parity.
The argument proceeds by induction on the number of squares of $\Delta$.
The construction of the quasi-perfect matching is trivial
if $\Delta$ has fewer than $4$ squares.

\begin{figure}[ht]
\psfrag{del}{$\delta$}
\psfrag{1}{$1$}
\psfrag{2}{$2$}
\psfrag{3}{$3$}
\begin{center}
\epsfig{height=30mm,file=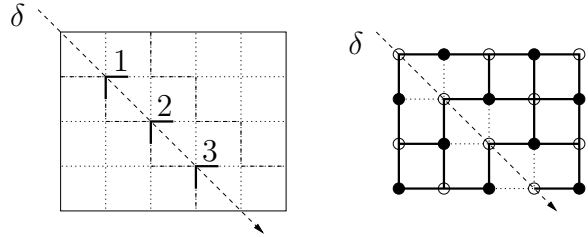}
\end{center}
\caption{Wedges along a good diagonal and the subgraph $\cG^\cR_\Delta$}
\label{fig:wedges}
\end{figure}

In general, start with a quadriculated disk $\Delta$ with $b=w$
and take a good diagonal $\delta$ as in Figure \ref{fig:wedges}.
Draw and number \textit{wedges} along $\delta$ as in the figure;
a tiling \textit{respects} a wedge if no domino in the tiling
crosses a leg of the wedge.
We define a partition $\cT_\Delta = \cD_\Delta \sqcup \cR_\Delta$:
a tiling $\tau$ belongs to $\cD_\Delta$
if and only if $\tau$ \textit{disrespects}
at least one of the wedges along $\delta$
(see \cite{Pachter} for a similar construction with a different purpose).
The loner of the quasi-perfect matching,
if it exists, will belong to $\cR_\Delta$;
the sets $\cD_\Delta$ and $\cR^\ast_\Delta = \cR_\Delta \cap \cT^\ast_\Delta$
will be invariant by $\rho$.
Equivalently, deletion of the edges of $\cG_\Delta$
crossing the wedges obtains a subgraph $\cG^\cR_\Delta$:
tilings in $\cT_\Delta$ (resp. $\cR_\Delta$)
correspond to matchings in $\cG_\Delta$ (resp. $\cG^\cR_\Delta$).

We first construct the restriction $\rho|_{\cD_\Delta}$.
Given $\tau \in \cD_\Delta$, assume that the first wedge to be disrespected
is the $k$-th wedge.
This means that the first $2\times 2$ square
formed by dominos along $\delta$ is positioned around that wedge:
$\rho(\tau)$ differs from $\tau$ by a flip in that square.

There is a natural bijection $\partial: \cR_\Delta \to \cT_{\Delta'}$,
where $\Delta'$ is the disjoint union of quadriculated disks
obtained from $\Delta$ by cut-and-paste along $\delta$. 
Indeed, for $\tau \in \cR_\Delta$,
define $\partial(\tau) \in \cT_{\Delta'}$
by removing the dominos covering one of the squares $s_{i+1/2}$
along $\delta$ and gluing the remaining parts.
Given a quasi-perfect matching
$\rho': \cT^\ast_{\Delta'} \to \cT^\ast_{\Delta'}$,
define $\cR^\ast_{\Delta} = \partial^{-1}(\cT^\ast_{\Delta'})$
and $\rho(\tau) = \partial^{-1}(\rho'(\partial(\tau)))$.

\begin{figure}[ht]
\psfrag{par}{$\partial$}
\psfrag{rho}{$\rho$}
\psfrag{rhop}{$\rho'$}
\psfrag{rhopp}{$\rho''$}
\begin{center}
\epsfig{height=60mm,file=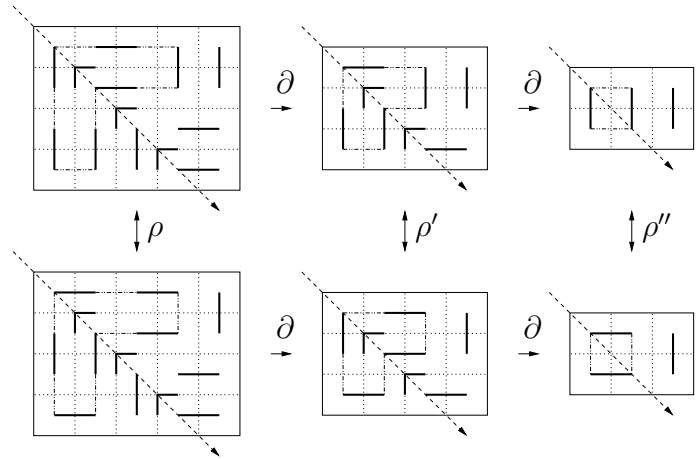}
\end{center}
\caption{The maps $\partial$ and $\rho$}
\label{fig:longcycle}
\end{figure}

If $\Delta'$ is a quadriculated disk,
a quasi-perfect matching is obtained by recursion.
Otherwise, for the detached collection
\[ \Delta' = \Delta'_1 \sqcup \cdots \sqcup \Delta'_d, \quad d > 1, \]
assume (again by recursion) that quasi-perfect matchings
$\rho'_i: \cT^\ast_{\Delta'_i} \to \cT^\ast_{\Delta'_i}$
have been obtained for each $\Delta'_i$ (possibly with loners).
For $\tau' \in \cT_{\Delta'}$,
let $\tau'_i$ be the restriction of $\tau'$ to $\Delta'_i$.
In order to find $\rho'(\tau')$,
search for the smallest $i$ for which $\tau'_i \in \cT^\ast_{\Delta'_i}$
(i.e., $\tau'_i$ is \textit{not} a loner);
construct $\rho'(\tau')$ by changing $\tau'$ in $\Delta'_i$ only:
\[ \rho'(\tau') = \tau'_1 \sqcup \cdots \sqcup \rho'_i(\tau'_i)
\sqcup \cdots \sqcup \tau'_d. \]
A tiling remains unmatched if and only if its restriction
to each $\Delta'_i$ is a loner:
since there is at most one loner in each $\cT_{\Delta'_i}$,
there is at most one loner in $\cT_{\Delta'}$
and $\rho'$ is indeed a quasi-perfect matching.

If the diagonal $\delta$ is unbalanced then $b' \ne w'$ and $\Delta'$
admits no domino tilings.
Consistently, in this case, $\cR_\Delta$ is empty:
this follows from the impossibility of respecting the last wedge.
More generally, if $\Delta' = \Delta'_1 \sqcup \cdots \sqcup \Delta'_d$
and (at least) one of the disks $\Delta'_i$ admits no domino tilings
then $\cR_\Delta$ is empty and we are done.

We must perform a final check:
$\tau$ and $\rho(\tau)$ are supposed to have opposite parities.
This is clear for $\tau \in \cD_\Delta$;
before we address the issue for $\tau \in \cR_\Delta$,
we present a few examples.

We follow the construction above in order to compute $\rho(\tau)$
where $\tau \in \cT_\Delta$ sits at the upper left hand corner
of Figure \ref{fig:longcycle}.
Recall that the definition of $\rho$ is dependent
on a specific choice of good diagonal not only for the original
disk $\Delta$ but for every disk reached in the process.
For $\delta$ as indicated, $\tau \in \cR_\Delta$.
Take $\tau' = \partial(\tau) \in \cT_{\Delta'}$
and a good diagonal $\delta'$ of $\Delta'$.
Again, $\tau' \in \cR_{\Delta'}$ so we must go to $\Delta''$
where $\tau'' = \partial(\tau') \in \cD_{\Delta''}$.
We construct $\rho''(\tau'') \in \cT_{\Delta''}$
(vertical arrow) and bring it back to obtain
$\rho'(\tau') = \partial^{-1}(\rho''(\tau'')) \in \cR_{\Delta'}$
and finally
$\rho(\tau) = \partial^{-1}(\rho'(\tau')) \in \cR_\Delta$.

\begin{figure}[ht]
\psfrag{par}{$\partial$}
\begin{center}
\epsfig{height=25mm,file=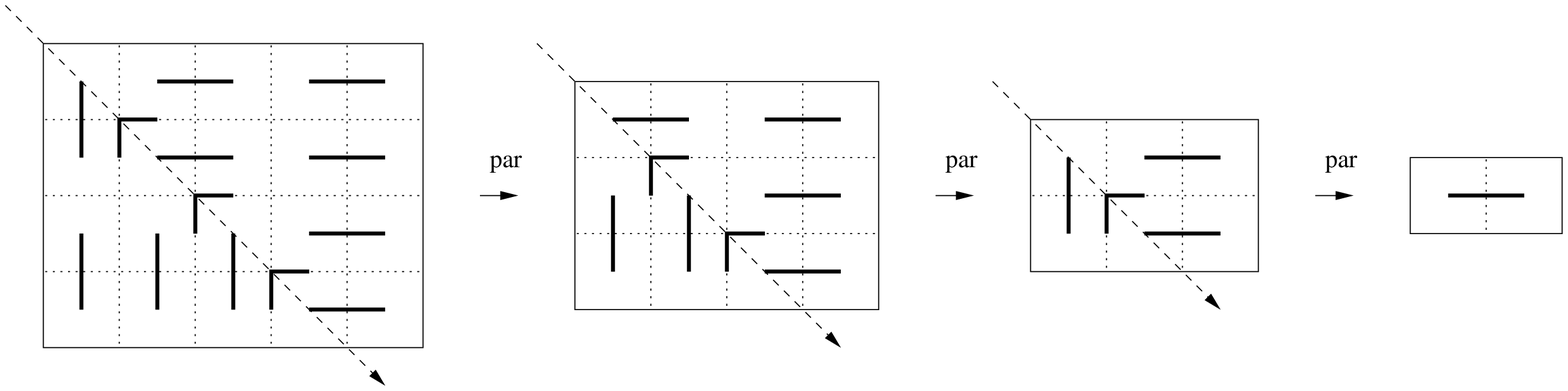}
\end{center}
\caption{A loner}
\label{fig:loner}
\end{figure}

\begin{figure}[ht]
\psfrag{par}{$\partial$}
\psfrag{rho}{$\rho$}
\psfrag{rhop}{$\rho'$}
\psfrag{rhopp}{$\rho''$}
\psfrag{rhoppp}{$\rho'''$}
\begin{center}
\epsfig{height=50mm,file=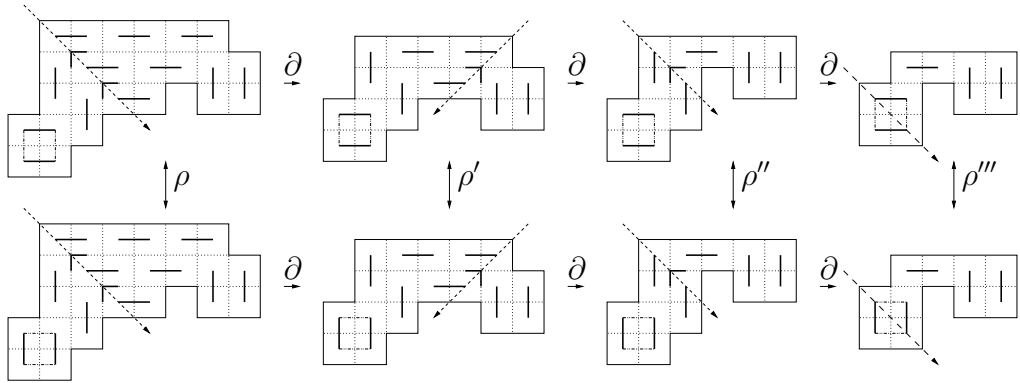}
\end{center}
\caption{Matching tilings in a more degenerate situation}
\label{fig:disk-w}
\end{figure}

In Figure \ref{fig:loner},
a loner is identified by a sequence of cut-and-paste operations
leading to a disk with a unique tiling.
In Figure \ref{fig:disk-w} we again compute $\rho(\tau)$
($\tau$ sits on the upper left corner);
notice that there is a large region where domino position is forced
but the construction still applies.

We recall some well known constructions.
The superposition $[\tau_1-\tau_2]$ of two tilings $\tau_1$ and $\tau_2$
consists of disjoint non-oriented simple closed curves
of consecutive dominos (or edges) alternating between $\tau_1$ and $\tau_2$;
dominos which are common to $\tau_1$ and $\tau_2$ are discarded.
Such curves are cycles (in a different sense) in the dual graph $\cG_\Delta$
but we reserve the word for permutation cycles.
Consider the bijections
$\pi_1, \pi_2: \{1, 2, \ldots, w \} \to \{1, 2, \ldots, b \}$
associated with the tilings $\tau_1, \tau_2$ and
decompose the permutation $\pi_2^{-1} \pi_1 \in S_w$
as a product of disjoint cycles.
These cycles correspond to the curves in $[\tau_1 - \tau_2]$
and the length of each curve (defined as the number of edges in $\cG_\Delta$)
is twice the length of the cycle.
The discarded dominos correspond to trivial cycles of length $1$ and
are irrelevant for parity checks.

If $\tau_1$ and $\tau_2$ differ by a flip then
$[\tau_1 - \tau_2]$ is a single curve of length $4$
and $\pi_2^{-1}\pi_1$ is a cycle of length $2$.
More generally,
two tilings $\tau_1, \tau_2 \in \cT_\Delta$ are \textit{compatible}
if $[\tau_1-\tau_2]$ consists of a single curve
whose length is a multiple of $4$;
we denote compatibility by $\tau_1 \leftrightarrow \tau_2$.
If $\tau_1 \leftrightarrow \tau_2$ then $\pi_2^{-1}\pi_1$
is a cycle of even length, an odd permutation,
and $\tau_1$ and $\tau_2$ have opposite parities.
We claim that, for $\tau_1, \tau_2 \in \cR_\Delta$,
\[ \tau_1 \leftrightarrow \tau_2 \quad\iff\quad
 \partial(\tau_1) \leftrightarrow \partial(\tau_2). \]
By the inductive construction of $\rho$,
the claim implies that $\tau_1 \leftrightarrow \rho(\tau_1)$,
completing the parity check.

\begin{figure}[ht]
\psfrag{par}{$\partial$}
\begin{center}
\epsfig{height=28mm,file=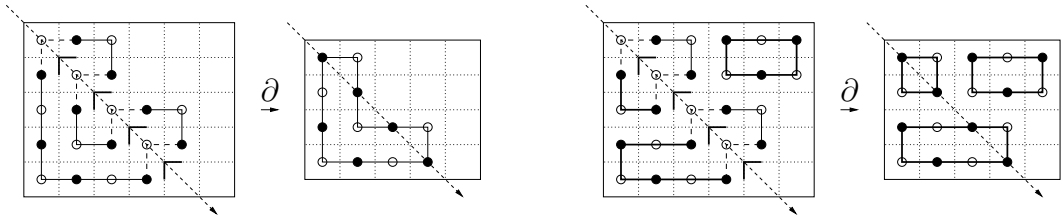}
\end{center}
\caption{Compatibility is preserved by cut-and-paste}
\label{fig:tunnel}
\end{figure}

Figure \ref{fig:tunnel} provides two examples of
$[\tau_1 - \tau_2]$ and $[\partial(\tau_1) - \partial(\tau_2)]$
for tilings $\tau_i \in \cR_\Delta$.
The reader should check that in the first example,
$\tau_1 \leftrightarrow \tau_2$ and 
$\partial(\tau_1) \leftrightarrow \partial(\tau_2)$;
in the second,
$\tau_1 \not\leftrightarrow \tau_2$ and 
$\partial(\tau_1) \not\leftrightarrow \partial(\tau_2)$.
Some vertices of the dual graphs $\cG_\Delta$ and $\cG_{\Delta'}$
are indicated for clarity.

In general, decompose the curves forming $[\tau_1 - \tau_2]$ into
dashed segments through \textit{corridors} between wedges
and solid segments on each side of the good diagonal.
Cut-and-paste deletes dashed segments 
and acts on solid segments by translation.
Thus, following solid segments yields a natural one-to-one correspondence
between curves in $[\tau_1 - \tau_2]$
and curves in $[\partial(\tau_1) - \partial(\tau_2)]$.
Furthermore, corresponding curves differ by the deletion
of dashed segments of length $2$, the passages of the curve through corridors.
Since at each such passage the curve goes
from one side of the diagonal to the other,
the number of passages for each curve is even.
Thus, lengths of corresponding curves are congruent $\mod 4$,
proving the claim and completing the proof.




\section{Algebraic cut-and-paste}
\label{sect:decomp}

The bulk of this section is dedicated to relating
the black-to-white matrices $B_\Delta$ and $B_{\Delta'}$
where $\Delta'$ is obtained from $\Delta$ by cut-and-paste
(there is no difficulty in defining black-to-white matrices for bicolored
disjoint union of quadriculated disks).
More precisely, assume that $\Delta$ (resp. $\Delta'$)
has $b$ (resp. $b'$) black squares and $w$ (resp. $w'$) white squares.
Let $I_n$ be the $n \times n$ identity matrix and
$I_{n,m}$ be the $n \times m$ defective identity matrix
with $(i,j)$ entry equal to 1 if $i=j$ and 0 otherwise.
We obtain in Lemma \ref{lemma:LDUstep} a factorization
\[ B_\Delta =
L_\Delta
\begin{pmatrix} I_{b-b',w-w'} & 0 \\ 0 & B_{\Delta'} \end{pmatrix}
U_\Delta  \]
where $L_\Delta$ and $U_\Delta$ are very special square triangular matrices.
This factorization is the inductive step in the proof of Theorem \ref{theo:LDU}.
We first present an example.

\begin{figure}[ht]
\begin{center}
\psfrag{del}{$\delta$}
\epsfig{height=45mm,file=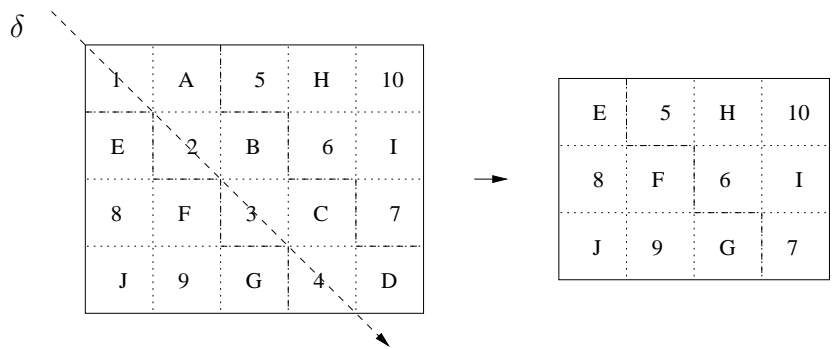}
\end{center}
\caption{Disks $\Delta$ and $\Delta'$}
\label{fig:rectangle}
\end{figure}

The quadriculated disks shown in Figure \ref{fig:rectangle}
have black-to-white matrices
\[ B_\Delta
= \begin{pmatrix}
1& & & &1& & & & & \\
1&1& & &1&1& & & & \\
 &1&1& & &1&1& & & \\
 & &1&1& & &1& & & \\
1&1& & & & & &1& & \\
 &1&1& & & & &1&1& \\
 & &1&1& & & & &1& \\
 & & & &1&1& & & &1\\
 & & & & &1&1& & &1\\
 & & & & & & &1&1&
\end{pmatrix}, \qquad
B_{\Delta'} = \begin{pmatrix}
1&1& &1& & \\
 &1&1&1&1& \\
 & &1& &1& \\
1&1& & & &1\\
 &1&1& & &1\\
 & & &1&1&
\end{pmatrix}. \]
Rows and columns are indexed by numbers and letters
respectively in Figure \ref{fig:rectangle}.
The first four rows and columns of $B_\Delta$
correspond to the eight squares removed by cut-and-paste.
Partition $B_\Delta$ in four blocks so that 
$B_{11} = B_\delta$ is the black-to-white matrix
of the disk around the diagonal $\delta$
consisting of squares $1, 2, 3, 4, A, B, C, D$ and
$B_{22}$ is the bottom $6\times 6$ principal minor.
Notice that $B_{22}$ and $B_{\Delta'}$ are very similar:
the difference lies in the top $3\times 3$ principal minor
of each matrix.
These positions describe adjacencies between squares
$5, 6, 7$ and $E, F, G$.

Elementary operations in rows and columns specified by
\[ \tilde X = \begin{pmatrix} 0&1&0&0 \\ 0&0&1&0 \\ 0&0&0&1 \\
0&0&0&0 \\ 0&0&0&0 \\ 0&0&0&0 \end{pmatrix}, \qquad
\tilde Y = \begin{pmatrix} 1&0&0&0&0&0 \\ 0&1&0&0&0&0 \\
0&0&1&0&0&0 \\ 0&0&0&0&0&0 \end{pmatrix} \]
can be applied to $B_\Delta$ to obtain a block diagonal matrix
\[ B_\Delta = \begin{pmatrix} I_4 & 0 \\ \tilde X & I_6 \end{pmatrix}
\begin{pmatrix} B_\delta & 0 \\ 0 & \tilde B_{\Delta'} \end{pmatrix}
\begin{pmatrix} I_4 & \tilde Y \\ 0 & I_6 \end{pmatrix}, \]
where
\[ \tilde B_{\Delta'} = \begin{pmatrix}
-1&-1& &1& & \\
&-1&-1&1&1& \\
& &-1& &1& \\
1&1& & & &1\\
&1&1& & &1\\
& & &1&1&
\end{pmatrix} \]
is surprisingly similar to $B_{\Delta'}$.
More precisely, $\tilde B_{\Delta'} = S_{b'} B_\Delta S_{w'}$
where $S_{b'} = \diag(-1,-1,-1,1,1,-1)$ and
$S_{w'} = \diag(1,1,1,-1,-1,1)$.
It is this ``coincidence'' that allows for this construction
to be used as the inductive step in the proof of Theorem \ref{theo:LDU}.

Before discussing the relationship between $B_\Delta$ and $B_{\Delta'}$
we present a lemma in linear algebra.
The proof is a straightforward computation left to the reader.

\begin{lemma}
\label{lemma:linalg}
Decompose an $(n + m) \times (n' + m')$ matrix $M$ as
\[ M = \begin{pmatrix} M_{11} & M_{12} \\ M_{21} & M_{22} \end{pmatrix}, \]
where $M_{11}$ is $n \times n'$.
If $n' \le n$ and $N$ is a $n' \times m'$ matrix with $M_{11} N = M_{12}$ then
\[ M =
\begin{pmatrix} M_{11} I_{n',n} & 0 \\ M_{21} I_{n',n} & I_{m} \end{pmatrix}
\begin{pmatrix} I_{n,n'} & 0 \\ 0 & M_{22} - M_{21} N \end{pmatrix}
\begin{pmatrix} I_{n'} & N \\ 0 & I_{m'} \end{pmatrix}.
\]
Similarly, if $n' \ge n$ and $N$ is a $m \times n$ matrix
with $N M_{11} = M_{21}$ then
\[ M =
\begin{pmatrix} I_{n} & 0 \\ N & I_{m} \end{pmatrix}
\begin{pmatrix} I_{n,n'} & 0 \\ 0 & M_{22} - N M_{12} \end{pmatrix}
\begin{pmatrix} I_{n',n} M_{11} &
I_{n',n} M_{12} \\ 0 & I_{m'} \end{pmatrix}. 
\]
\end{lemma}

The next lemma is the inductive step in the proof of Theorem \ref{theo:LDU}.

\begin{lemma}
\label{lemma:LDUstep}
Let $\Delta$ be a quadriculated disk with $b$ black and $w$ white squares,
$b+w > 1$.
Let $\Delta' = \Delta'_1 \sqcup \cdots \sqcup \Delta'_d$
(with $b' = b'_1 + \cdots + b'_d$ black
and $w' = w'_1 + \cdots + w'_d$ white squares)
be obtained from $\Delta$ by cut-and-paste along a good diagonal $\delta$.
Label black and white squares in $\Delta$ so that removed squares
come first, in the order prescribed by the good diagonal;
label squares in $\Delta'$ next.
Then the black-to-white matrices $B_\Delta$ and $B_{\Delta'}$ satisfy
\[ B_\Delta =
\begin{pmatrix} L & 0 \\ X & S_{b'} \end{pmatrix}
\begin{pmatrix} I_{b-b',w-w'} & 0 \\ 0 & B_{\Delta'} \end{pmatrix}
\begin{pmatrix} U & Y \\ 0 & S_{w'} \end{pmatrix}
\]
where $L$ (resp. $U$) is an invertible lower (resp. upper)
square matrix of order $b-b'$ (resp. $w-w'$)
and $S_{b'}$ and $S_{w'}$ are square diagonal matrices.
Furthermore, all entries of $S_{b'}$, $S_{w'}$, $L$, $U$, $X$ and $Y$
equal $0$, $1$ or $-1$.
\end{lemma}

The statement above requires clarification in some degenerate cases.
If $\Delta'$ is empty, $B_{\Delta'}$ collapses and
$B_\Delta = L I_{b,w} U$.
If instead $\Delta'$ is a disjoint union of unit squares,
all of the same color, then either $w' = 0$ or $b' = 0$ and
\[ B_\Delta =
\begin{pmatrix} L & 0 \\ X & S_{b'} \end{pmatrix}
\begin{pmatrix} I_{b-b',w} \\ 0 \end{pmatrix} U
\quad\hbox{or}\quad
B_\Delta =
L \begin{pmatrix} I_{b,w-w'} & 0 \end{pmatrix}
\begin{pmatrix} U & Y \\ 0 & S_{w'} \end{pmatrix}.
\]

{\nobf Proof:}
Assume that the deleted squares are $s_{1/2}, s^l_1, \ldots$
and that the square $s_{1/2}$ is black; thus $k = b-b'$;
if $s_{1/2}$ were white all computations would be transposed.
Let $j_1, \ldots, j_{k-1}$ be the indices of the white squares
$s^r_1, \ldots, s^r_{k-1}$; notice that $j_i > w-w'$.
Decompose the matrix $B_\Delta$ in four blocks,
\[ B_\Delta = \begin{pmatrix} B_{11} & B_{12} \\ B_{21} & B_{22} \end{pmatrix}, \]
where $B_{22}$ is a $b' \times w'$ matrix.
By construction, $B_{11}$ has one of the two forms below,
the first case corresponding to balanced good diagonals
(i.e., to $b-b' = w-w'$).
\[
B_{11} = \begin{pmatrix}
1 & 0 & 0 & \cdots & 0 & 0 \\
1 & 1 & 0 & \cdots & 0 & 0 \\
0 & 1 & 1 & \cdots & 0 & 0 \\
\vdots & \vdots & \vdots & & \vdots & \vdots \\
0 & 0 & 0 & \cdots & 1 & 0 \\
0 & 0 & 0 & \cdots & 1 & 1
\end{pmatrix},
\quad \hbox{or} \quad
B_{11} = \begin{pmatrix}
1 & 0 & 0 & \cdots & 0  \\
1 & 1 & 0 & \cdots & 0  \\
0 & 1 & 1 & \cdots & 0  \\
\vdots & \vdots & \vdots & & \vdots \\
0 & 0 & 0 & \cdots & 1  \\
0 & 0 & 0 & \cdots & 1 
\end{pmatrix}.
\]
Let $S_b$ (resp. $S_w$) be a $b \times b$ (resp. $w \times w$)
diagonal matrix with diagonal entries equal to $1$ or $-1$;
the $i$-th entry of $S_b$ (resp. $S_w$) is $-1$ if the $i$-th
black (resp. white) square is strictly to the right of $\delta$.
Write
\[ S_b = \begin{pmatrix} I_{b-b'} & 0 \\ 0 & S_{b'} \end{pmatrix}, \qquad
S_w = \begin{pmatrix} I_{w-w'} & 0 \\ 0 & S_{w'} \end{pmatrix}. \]
We have
\[ S_b B_\Delta S_w =
\begin{pmatrix} B_{11} & - B_{12} \\ B_{21} & B_{22} \end{pmatrix}. \]
The nonzero entries of $B_{12}$ are
$(i,j_i)$ and $(i+1,j_i)$ for $i = 1,\ldots,k-1$.
Thus, the nonzero columns of $B_{12}$ equal to the first $k-1$ columns
of $B_{11}$. Let $N$ be the $(w-w') \times w'$ matrix with entries $0$ or $-1$,
with nonzero entries at $(1,j_1), (2,j_2), \ldots, (k-1,j_{k-1})$.
Clearly $B_{11} N = - B_{12}$ and we may apply Lemma \ref{lemma:linalg}
to write
\[ S_b B_\Delta S_w =
\begin{pmatrix} B_{11} I_{w-w',b-b'} & 0 \\
B_{21} I_{w-w',b-b'} & I_{b'} \end{pmatrix}
\begin{pmatrix} I_{b-b',w-w'} & 0 \\ 0 & B_{22} - B_{21} N \end{pmatrix}
\begin{pmatrix} I_{w-w'} & N \\ 0 & I_{w'} \end{pmatrix}.
\]
We claim that $B_{\Delta'} = B_{22} - B_{21} N$.
The nonzero columns of the matrix $- B_{21} N$ are
the columns of $B_{21}$, except that the first column is moved
to position $j_1$, the second column is moved to $j_2$ and so on.
The $k$-th column of $B_{21}$, if it exists, is discarded.
These nonzero entries correspond precisely to the identifications
which must be performed in order to obtain $\Delta'$, i.e.,
to the ones which must be added to $B_{22}$ in order to obtain $B_{\Delta'}$.
Clearing up signs,
\[ B_\Delta =
\begin{pmatrix} B_{11} I_{w-w',b-b'} & 0 \\
S_{b'} B_{21} I_{w-w',b-b'} & S_{b'} \end{pmatrix}
\begin{pmatrix} I_{b-b',w-w'} & 0 \\ 0 & B_{\Delta'} \end{pmatrix}
\begin{pmatrix} I_{w-w'} & N S_{w'} \\ 0 & S_{w'} \end{pmatrix}.
\]
If the good diagonal is balanced, this finishes the proof.
In the unbalanced case, $\tilde L = B_{11} I_{w-w',b-b'}$ is not invertible
since its last column is zero. Replace the $(k,k)$ entry of 
$\tilde L$ by 1 to obtain a new matrix $L$:
$L$ is clearly invertible and $\tilde L I_{b-b',w-w'} = L I_{b-b',w-w'}$.
The proof is now complete.
\qed

{\nobf Proof of Theorem \ref{theo:LDU}:}
The basis of the induction on the number of squares of $\Delta$
consists of checking that the theorem holds for disks with at most two squares.
Notice that if the disk consists of a single square then $b = 0$ or $w = 0$
and the matrices are degenerate.

Let $\Delta$ be a quadriculated disk and
$\Delta' = \Delta'_1 \sqcup \cdots \sqcup \Delta'_d$ 
be obtained from $\Delta$ by cut-and-paste.
By induction on the number of squares the theorem may be assumed to hold
for eack $\Delta'_k$ and we therefore write
$B_{\Delta'} = L_{\Delta'} \tD_{\Delta'} U_{\Delta'}$.
From the induction step, Lemma \ref{lemma:LDUstep}, write
\begin{align}
B_\Delta &=
\begin{pmatrix} L_{\step} & 0 \\ X_{\step} & S_{b'} \end{pmatrix}
\begin{pmatrix} I_{b-b',w-w'} & 0 \\ 0 & B_{\Delta'} \end{pmatrix}
\begin{pmatrix} U_{\step} & Y_{\step} \\ 0 & S_{w'} \end{pmatrix}
\notag\\
&= L_\Delta \tD_\Delta U_\Delta. \notag
\end{align}
where
\[
L_\Delta =
\begin{pmatrix} L_{\step} & 0 \\ X_{\step} & S_{b'} \end{pmatrix}
\begin{pmatrix} I_{b-b'} & 0 \\ 0 & L_{\Delta'} \end{pmatrix},
\quad
U_\Delta =
\begin{pmatrix} I_{w-w'} & 0 \\ 0 & U_{\Delta'} \end{pmatrix}
\begin{pmatrix} U_{\step} & Y_{\step} \\ 0 & S_{w'} \end{pmatrix}.
\]
The theorem now follows from observing that each nonzero entry of $L_\Delta$
(resp. $U_\Delta$)
is, up to sign, copied from either
$L_{\Delta'}$, $L_{\step}$ or $X_{\step}$
(resp. $U_{\Delta'}$, $U_{\step}$ or $Y_{\step}$)
and is therefore equal to $1$ or $-1$.
\qed

We present a direct consequence of Theorem \ref{theo:LDU}.

\begin{coro}
\label{coro:tor}
Let $\Delta$ be a quadriculated disk with black-to-white matrix $B_\Delta$.
If $v$ has integer entries and the system $B_\Delta x = v$ admits
a rational solution then the system admits an integer solution.
\end{coro}

This corollary may be interpreted as saying that the co-kernel
$\ZZ^b/B_\Delta(\ZZ^w)$ of $B_\Delta: \ZZ^w \to \ZZ^b$
is a free abelian group.
From Theorem \ref{theo:LDU},
the rank $r$ of $B_\Delta$ is the same in $\QQ$ as in $\ZZ_p$
for any prime number $p$.
Notice that the proof of Theorem \ref{theo:DT} in \cite{DT}
is based on this fact for $p=2$.

\begin{figure}[ht]
\begin{center}
\epsfig{height=30mm,file=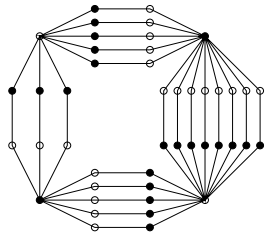}
\end{center}
\caption{Determinant 1 does not imply $L\tilde DU$ decomposition}
\label{fig:fibo}
\end{figure}

The example in Figure \ref{fig:fibo} is instructive:
the $B_G$ matrix of this planar graph $G$ has determinant 1
but admits no $L\tilde DU$ decomposition where the matrices
have integer coefficients since the removal of any two vertices
of opposite colors from $G$ yields a graph whose determinant
has absolute value greater than 1.

\section{Boards}
\label{sect:boards}

\begin{figure}[ht]
\begin{center}
\epsfig{height=30mm,file=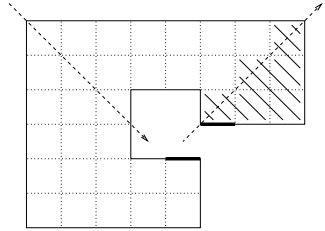}
\end{center}
\caption{A board and two good diagonals, one excellent.}
\label{fig:board}
\end{figure}

Topological subdisks of $\RR^2$
consisting of unit squares with vertices in $\ZZ^2$ are {\it boards}.
In other words, a board is a topological subdisk of $\RR^2$
whose boundary is a polygonal curve consisting of segments
of length 1 joining points in $\ZZ^2$.
The quadriculated disk in Figure \ref{fig:disk0} is not a board.
The class of boards is not closed under cut-and-paste:
in Figure \ref{fig:board}, the two enhanced segments on the boundary
would be superimposed by cut-and-paste along the good diagonal on the left.
Cut-and-paste along the good diagonal indicated on the right, however,
yields a smaller board.
The main result of this section is that, given a board $\Delta$,
it is always possible to choose a good diagonal $\delta$
such that $\Delta' = \partial_\delta(\Delta)$ is a disjoint union of boards.

\begin{figure}[ht]
\psfrag{Positive}{Positive}
\psfrag{Negative}{Negative}
\begin{center}
\epsfig{height=15mm,file=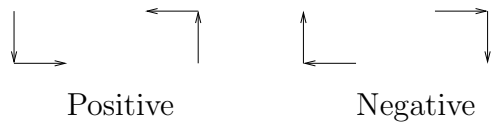}
\end{center}
\caption{Positive and negative boundary vertices}
\label{fig:posneg}
\end{figure}

Orient the boundary of a board $\Delta$ counterclockwise,
so that $\Delta$ lies to the left of the boundary.
Consider boundary vertices which are local extrema for
the restriction of $x+y$ to the boundary:
as in Figure \ref{fig:posneg},
call such vertices \textit{positive} if they are corners
(equivalently, if they are local extrema
for the restriction of $x+y$ to $\Delta$)
and \textit{negative} otherwise.
Let $V_{B,+}$ (resp. $V_{B,-}$) be the number of positive (resp. negative)
boundary vertices.

\begin{lemma}
\label{lemma:gaussbonnetne}
$V_{B,+} - V_{B,-} = 2$.
\end{lemma}

{\nobf Proof:}
Define $F$, $E$, $E_I$, $E_B$ and $V_I$ as in Lemma \ref{lemma:gaussbonnet}.
The number of boundary vertices is $V_B = V_{B,+} + V_{B,-} + V_{B,0}$
where $V_{B,0}$ is the number of boundary vertices
which are neither positive nor negative.
For each square, consider its $NW$ and $SE$ vertices:
interior vertices and negative vertices are counted twice,
positive vertices are not counted and
other boundary vertices are counted once and therefore
$2F = 2V_I + 2V_{B,-} + V_{B,0} = 2V_I + V_B - (V_{B,+} - V_{B,-})$.
Recall that $4F = 2E - E_B$ (Lemma \ref{lemma:gaussbonnet})
and $E = V + F - 1$ (Euler) and therefore
$2F = 2V_I + V_B - 2$, completing the proof.
\qed

\begin{theo}
\label{theo:board}
It is always possible to cut-and-paste a given board $\Delta$
to obtain a disjoint union of boards $\Delta'$.
\end{theo}

{\nobf Proof:}
A diagonal is \textit{excellent} if the $x$ and $y$ coordinates
are both monotonic along one of the two boundary arcs between $v_0$ and $v_k$;
without loss, let this arc lie to the right of the diagonal.
Excellent diagonals are good: the vertex $v^r_{k-1/2}$ is on the boundary.
We interpret cut-and-paste along an excellent diagonal
as leaving $\Delta^l$ fixed and moving $\Delta^r$.
In this way, $\Delta'$ becomes a subset of $\Delta$
and is therefore a disjoint union of boards.
We are left with proving that any board admits excellent diagonals.
Each diagonal defines two boundary arcs:
order these arcs by inclusion.
We claim that a diagonal defining a minimal arc is excellent.

Let $\delta^m = (v^m_0v^m_1\ldots v^m_k)$
be a diagonal inducing a minimal arc $\alpha$:
assume without loss of generality that $v^m_i = (a+i,b+i)$ 
for integers $a$ and $b$.
Consider the set $\tilde\Delta$ (dashed in Figure \ref{fig:board})
consisting of the squares totally or partially surrounded
by $\alpha$ and $\delta^m$.
It is easy to verify that $\tilde\Delta$ is a legitimate board
with boundary consisting of $\alpha$ and $\zeta$,
where $\zeta$ is the zig-zag line next to $\delta^m$.
Thus, the last edge of $\zeta$ can not overlap with $\alpha$
without contradicting the fact that the boundary point $v^m_k$
of $\Delta$ is surrounded by at most three squares in $\Delta$.

By Lemma \ref{lemma:gaussbonnetne},
the board $\tilde\Delta$ has at least two positive boundary points.
We claim that the existence of a positive boundary point
distinct from $v^m_0$ and $v^m_k$ contradicts minimality.
Notice that at this point it is clear that $v^m_0$ is positive;
the status of $v^m_k$ as a positive boundary point
will only follow from the claim.
Indeed, such a positive point $\hat v$
can not belong to the zig-zag line $\zeta$
and must therefore belong to $\alpha$.
Draw a diagonal $\hat\delta$ starting at $\hat v$:
being parallel to $\zeta$, $\hat\delta$ must intersect
the boundary of $\tilde\Delta$ in $\alpha$
and therefore defines a smaller arc $\hat\alpha$,
contradicting minimality and proving the claim.
Again by Lemma \ref{lemma:gaussbonnetne},
there are no negative boundary vertices.
In particular, there are
no positive or negative boundary vertices in $\alpha$
and we are done.
\qed

%



\bigbreak
\bigbreak

\bigskip\bigskip\bigbreak

{

\parindent=0pt
\parskip=0pt
\obeylines

Nicolau C. Saldanha and Carlos Tomei 
Departamento de Matem\'atica, PUC-Rio
R. Marqu\^es de S. Vicente 225, Rio de Janeiro, RJ 22453-900, Brazil

\medskip

nicolau@mat.puc-rio.br; http://www.mat.puc-rio.br/$\sim$nicolau/
tomei@mat.puc-rio.br

}

\end{document}